\documentclass[10pt,a4paper]{amsart}

\usepackage{amssymb}
\usepackage{amsthm}
\usepackage{amsfonts}
\usepackage{microtype}
\usepackage{mathrsfs}
\usepackage{graphicx}
\usepackage{color}
\usepackage{latexsym}
\usepackage{rotating}
\usepackage{xspace}
\usepackage[all]{xy}
\usepackage{longtable}
\usepackage{leftidx}
\usepackage{mathtools}
\usepackage{titlesec}
\usepackage{tikz}
\usetikzlibrary{calc}
\usetikzlibrary{arrows}
\usetikzlibrary{decorations.markings}
\usepackage{geometry}

\newtheorem{theorem}{Theorem}[section]
\numberwithin{equation}{theorem}
\newtheorem{lemma}[theorem]{Lemma}

\newtheorem{corollary}[theorem]{Corollary}

\theoremstyle{definition}
\newtheorem{definition}[theorem]{Definition}
\newtheorem{example}[theorem]{Example}

\theoremstyle{conjecture}

\newtheorem{question}[theorem]{Question}
\newtheorem{acknowledgement}{Acknowledgement}

\newcommand{\im}{\operatorname{im}}

\newcommand{\ara}{\operatorname{ara}}

\newcommand{\amp}{\operatorname{amp}}

\newcommand{\cd}{\operatorname{cd}}

\newcommand{\id}{\operatorname{id}}

\newcommand{\V}{\operatorname{V}}

\newcommand{\Cone}{\operatorname{Cone}}

\newcommand{\Ext}{\operatorname{Ext}}
\newcommand{\Supp}{\operatorname{Supp}}

\newcommand{\Hom}{\operatorname{Hom}}

\newcommand{\coker}{\operatorname{coker}}

\newcommand{\lo}{\longrightarrow}

\newcommand{\fa}{\frak{a}}

\newcommand{\suchthat}{\;\ifnum\currentgrouptype=16 \middle\fi|\;}

\newenvironment{prf}[1][Proof]{\begin{proof}[\bf #1]}{\end{proof}}

\makeatletter
\newcommand{\holim@}[2]{%
  \vtop{\m@th\ialign{##\cr
    \hfil$#1\operator@font holim$\hfil\cr
    \noalign{\nointerlineskip\kern1.5\ex@}#2\cr
    \noalign{\nointerlineskip\kern-\ex@}\cr}}%
}
\newcommand{\holim}{%
  \mathop{\mathpalette\holim@{\rightarrowfill@\textstyle}}\nmlimits@
}
\makeatother

\makeatletter
\def\@secnumfont{\bfseries}
\def\section{\@startsection{section}{1}%
  \z@{.7\linespacing\@plus\linespacing}{.5\linespacing}%
  {\normalfont\Large\bfseries\filcenter}}
\def\subsection{\@startsection{subsection}{2}%
  \z@{.5\linespacing\@plus.7\linespacing}{-.5em}%
  {\normalfont\large\bfseries}}
\makeatother

\begin{document}

\author[K. Divaani-Aazar, H. Faridian and M. Tousi]{Kamran Divaani-Aazar, Hossein Faridian
and Massoud Tousi}

\title[A New Outlook on Cofiniteness]
{A New Outlook on Cofiniteness}

\address{K. Divaani-Aazar, Department of Mathematics, Alzahra University, Vanak, Post Code 19834, Tehran, Iran-and-School of Mathematics, Institute for
Research in Fundamental Sciences (IPM), P.O. Box 19395-5746, Tehran, Iran.}
\email{kdivaani@ipm.ir}

\address{H. Faridian, Department of Mathematics, Shahid Beheshti University, G.C., Evin, Tehran, Iran, Zip Code 1983963113.}
\email{h.faridian@yahoo.com}

\address{M. Tousi, Department of Mathematics, Shahid Beheshti University, G.C., Evin, Tehran, Iran, P.O. Box 19395-5746.}
\email{mtousi@ipm.ir}


\keywords {Abelian category; cofinite module; cohomological dimension; derived category; dualizing complex; local cohomology module; local homology module.\\
The research of the first author is supported by a grant from IPM (No. 95130212).}

\begin{abstract}
Let $\mathfrak{a}$ be an ideal of a commutative noetherian (not necessarily local) ring $R$. In the case $\cd(\mathfrak{a},R)\leq 1$, we show that the subcategory
of $\mathfrak{a}$-cofinite $R$-modules is abelian. Using this and the technique of way-out functors, we show that if $\cd(\mathfrak{a},R)\leq 1$, or $\dim(R/\mathfrak{a})
\leq 1$, or $\dim(R) \leq 2$, then the local cohomology module $H^{i}_{\mathfrak{a}}(X)$ is $\mathfrak{a}$-cofinite for every $R$-complex $X$ with finitely generated
homology modules and every $i \in \mathbb{Z}$. We further answer Question 1.3 in the three aforementioned cases, and reveal a correlation between Questions 1.1,
1.2, and 1.3.

\end{abstract}

\maketitle

\section{Introduction}

\sloppy

Throughout this paper, $R$ denotes a commutative noetherian ring with identity and $\mathcal{M}(R)$ flags the category of $R$-modules.

In 1969, Hartshorne introduced the notion of cofiniteness for modules and complexes; see \cite{Ha1}. He defined an $R$-module $M$ to be $\mathfrak{a}$-cofinite if $\Supp_{R}(M)\subseteq \V(\mathfrak{a})$ and $\Ext^{i}_{R}(R/\mathfrak{a},M)$ is finitely generated for every $i\geq 0$. Moreover, in the case where $R$ is an $\mathfrak{a}$-adically complete regular ring of finite Krull dimension, he defined an $R$-complex $X$ to be $\mathfrak{a}$-cofinite if $X \simeq {\bf R}\Hom_{R}\left(Y,{\bf R} \Gamma_{\mathfrak{a}}(R)\right)$ for some $R$-complex $Y$ with finitely generated homology modules. He then proceeded to pose three questions in this direction which we paraphrase as follows.

\begin{question} \label{1.1}
Is the local cohomology module $H^{i}_{\mathfrak{a}}(M)$, $\mathfrak{a}$-cofinite for every finitely generated $R$-module $M$ and every $i \geq 0$?
\end{question}

\begin{question} \label{1.2}
Is the category $\mathcal{M}(R,\mathfrak{a})_{cof}$ consisting of $\mathfrak{a}$-cofinite $R$-modules an abelian subcategory of $\mathcal{M}(R)$?
\end{question}

\begin{question} \label{1.3}
Is it true that an $R$-complex $X$ is $\mathfrak{a}$-cofinite if and only if the homology module $H_{i}(X)$ is $\mathfrak{a}$-cofinite for every $i\in \mathbb{Z}$?
\end{question}

By providing a counterexample, Hartshorne showed that the answers to these questions are negative in general; see \cite[Section 3]{Ha1}. However, he established affirmative answers to these questions in the case where $\mathfrak{a}$ is a principal ideal generated by a nonzerodivisor and $R$ is an $\mathfrak{a}$-adically complete regular ring of finite Krull dimension, and also in the case where $\mathfrak{a}$ is a prime ideal with $\dim(R/\mathfrak{a})=1$ and $R$ is a complete regular local ring; see \cite[Propositions 6.1 and 6.2, Corollary 6.3, Theorem 7.5, Proposition 7.6 and Corollary 7.7]{Ha1}. Since then many papers are devoted to study his first two questions; see for example \cite{HK}, \cite{DM}, \cite{Ka1}, \cite{Ka2}, \cite{Me1}, \cite{Me2} and \cite{Y}. These results were extended in several stages to take the following form:

\begin{theorem} \label{1.4}
Let $\mathfrak{a}$ be an ideal of $R$ such that either $\ara(\mathfrak{a}) \leq 1$, or $\dim(R/\mathfrak{a}) \leq 1$, or $\dim(R) \leq 2$. Then $H^{i}_{\mathfrak{a}}(M)$ is $\mathfrak{a}$-cofinite for every finitely generated $R$-module $M$ and every $i \geq 0$, and $\mathcal{M}(R,\mathfrak{a})_{cof}$ is an abelian subcategory of $\mathcal{M}(R)$.
\end{theorem}

For the case $\ara(\mathfrak{a}) \leq 1$, refer to \cite[Theorem 1]{Ka2} and \cite[Theorem 2.1]{Ka1}. For the case $\dim(R/\mathfrak{a})
\leq 1$, see \cite[Theorem 2.6 and Corollary 2.12]{Me1}, \cite[Corollary 2.8]{BNS}, and \cite[Corollary 2.7]{BN}. For the case $\dim(R) \leq 2$, observe \cite[Theorem 7.10]{Me2} and \cite[Theorem 7.4]{Me2}.

The significance of cofiniteness of the local cohomology modules mainly stems from the fact that if an $R$-module $M$ is $\frak{a}$-cofinite, then its set of associated primes
is finite as well as all its Bass numbers and Betti numbers with respect to every prime ideal of $R$. It is worth mentioning that the investigation of such finiteness properties is a long-sought problem in commutative and homological algebra; see e.g. \cite{HS} and \cite{Ly}.

In this paper, we deal with the above three questions. Theorems \ref{2.1}, \ref{3.3}, \ref{3.5} and \ref{4.3} are our main results.

In \cite[Question 1]{PAB}, the authors asked: Is $\mathcal{M}(R,\mathfrak{a})_{cof}$ an abelian subcategory of $\mathcal{M}(R)$ for every ideal $\frak{a}$ of $R$ with $\cd(\frak{a},R)\leq 1$? We answer this question affirmatively by deploying the theory of local homology; see Theorem \ref{2.1}. Note that there exists an inequality $\cd(\mathfrak{a},R)\leq \ara(\mathfrak{a})$
that can be strict; see Example \ref{2.2}.

It turns out that to establish the cofiniteness of $H^{i}_{\mathfrak{a}}(X)$ for any $R$-complex $X$ with finitely generated homology modules, all we need to know is the cofiniteness of $H^{i}_{\mathfrak{a}}(M)$ for any finitely generated $R$-module $M$ and the abelianness of $\mathcal{M}(R,\mathfrak{a})_{cof}$; see Theorem \ref{3.3}. The crucial step to achieve this is to recruit the technique of way-out functors.

To be consistent in both module and complex cases, we define an $R$-complex $X$ to be $\mathfrak{a}$-cofinite if $\Supp_{R}(X)\subseteq \V(\mathfrak{a})$ and ${\bf R}\Hom_{R}(R/\mathfrak{a},X)$ has finitely generated homology modules. Corollary \ref{4.2} indicates that, for homologically bounded $R$-complexes, this definition coincides with that of Hartshorne.

Questions \ref{1.1} and \ref{1.2} have been high-profile among researchers, whereas not much attention has been brought to Question \ref{1.3}. The most striking result on this question is \cite[Theorem 1]{EK} which confines itself to complete Gorenstein local domains and the case $\dim (R/\frak{a})=1$. We answer Hartshorne's third question in the cases $\cd(\mathfrak{a},R) \leq 1$, $\dim(R/\mathfrak{a})\leq 1$, and $\dim(R) \leq 2$ with no extra assumptions on $R$; see Corollary \ref{3.6} (ii). Having the results thus far obtained at our disposal, we show that the answers to Questions \ref{1.1} and \ref{1.2} are affirmative if and only if the answer to Question \ref{1.3} is affirmative for all homologically bounded $R$-complexes; see Theorem \ref{4.3}.

\section{Question 1.2}

We need to work in the framework of the derived category $\mathcal{D}(R)$. For more information, refer to \cite{AF}, \cite{Ha2}, \cite{Fo}, \cite{Li}, and \cite{Sp}.

We let $\mathcal{D}_{\sqsubset}(R)$ (res. $\mathcal{D}_{\sqsupset}(R)$) denote the full subcategory of $\mathcal{D}(R)$ consisting of $R$-complexes $X$ with $H_{i}(X)=0$
for $i \gg 0$ (res. $i \ll 0$), and let  $D_{\square}(R):=\mathcal{D}_{\sqsubset}(R)\cap \mathcal{D}_{\sqsupset}(R)$. We further let $\mathcal{D}^{f}(R)$ denote the full
subcategory of $\mathcal{D}(R)$ consisting of $R$-complexes $X$ with finitely generated homology modules. We also feel free to use any combination of the subscripts and
the superscript as in $\mathcal{D}^{f}_{\square}(R)$, with the obvious meaning of the intersection of the two subcategories involved.

\begin{lemma} \label{2.0}
Let $\mathfrak{a}$ be an ideal of $R$ and $X\in D_{\square}(R)$. Then the following conditions are equivalent:
\begin{enumerate}
\item[(i)] ${\bf R}\Hom_{R}(R/\mathfrak{a},X)\in D^f(R)$.
\item[(ii)] ${\bf L}\Lambda^{\mathfrak{a}}(X) \in \mathcal{D}^{f}_{\square}\left(\widehat{R}^{\mathfrak{a}}\right)$.
\end{enumerate}
\end{lemma}

\begin{prf}
See \cite[Propositions 7.4]{WW}.
\end{prf}

In this section, we show that given an ideal $\mathfrak{a}$ of $R$ with $\cd(\frak{a},R)\leq 1$, the subcategory $\mathcal{M}(R,\mathfrak{a})_{cof}$ of $\mathcal{M}(R)$
is abelian. This fact is proved in \cite[Theorem 2.4]{PAB}, under the extra assumption that $R$ is local. Here we relax this assumption. The tool here is the local homology
functors.

Recall that the local homology functors are the left derived functors of the completion functor. More precisely, $H^{\mathfrak{a}}_{i}(-):= L_{i}\left(\Lambda^{\mathfrak{a}}(-)
\right)$ for every $i \geq 0$, where $\Lambda^{\mathfrak{a}}(M):=\widehat{M}^{\mathfrak{a}}=\underset{n}\varprojlim (M/\mathfrak{a}^{n}M)$ for any $R$-module $M$. Further, we
remind the cohomological dimension of $M$ with respect to $\mathfrak{a}$ as $$\cd(\mathfrak{a},M):=\sup \left\{i \in \mathbb{Z} \suchthat H^{i}_{\mathfrak{a}}(M)\neq 0
\right\}.$$

\begin{theorem} \label{2.1}
Let $\mathfrak{a}$ be an ideal of $R$. Then the following assertions hold:
\begin{enumerate}
\item[(i)] An $R$-module $M$ with $\Supp_{R}(M)\subseteq V(\mathfrak{a})$ is $\mathfrak{a}$-cofinite if and only if $H^{\mathfrak{a}}_{i}(M)$ is a finitely generated $\widehat{R}^{\mathfrak{a}}$-module for every $0 \leq i \leq \cd(\mathfrak{a},R)$.
\item[(ii)] If $\cd(\mathfrak{a},R)\leq 1$, then $\mathcal{M}(R,\mathfrak{a})_{cof}$ is an abelian subcategory of $\mathcal{M}(R)$.
\end{enumerate}
\end{theorem}

\begin{prf}
(i): By \cite[Theorem 2.5 and Corollary 3.2]{GM}, $H^{\mathfrak{a}}_{i}(M)=0$ for every $i> \cd(\mathfrak{a},R)$. Therefore, the assertion follows from Lemma \ref{2.0}.

(ii): Let $M$ and $N$ be two $\mathfrak{a}$-cofinite $R$-modules and $f:M\rightarrow N$ an $R$-homomorphism. The short exact sequence
\begin{equation} \label{eq:2.1.1}
0 \rightarrow \ker f \rightarrow M \rightarrow \im f \rightarrow 0,
\end{equation}
gives the exact sequence
$$H^{\mathfrak{a}}_{0}(M) \rightarrow H^{\mathfrak{a}}_{0}\left(\im f \right) \rightarrow 0,$$
which in turn implies that $H^{\mathfrak{a}}_{0}\left(\im f \right)$ is finitely generated $\widehat{R}^{\mathfrak{a}}$-module
since $H^{\mathfrak{a}}_{0}(M)$ is so. The short exact sequence
\begin{equation} \label{eq:2.1.2}
0 \rightarrow \im f \rightarrow N \rightarrow \coker f \rightarrow 0,
\end{equation}
gives the exact sequence
\begin{equation} \label{eq:2.1.3}
H^{\mathfrak{a}}_{1}(N) \rightarrow H^{\mathfrak{a}}_{1}\left(\coker f \right) \rightarrow H^{\mathfrak{a}}_{0}\left(\im f \right)
\rightarrow H^{\mathfrak{a}}_{0}(N) \rightarrow H^{\mathfrak{a}}_{0}\left(\coker f \right) \rightarrow 0.
\end{equation}
As $H^{\mathfrak{a}}_{0}(N)$, $H^{\mathfrak{a}}_{0}\left(\im f \right)$, and $H^{\mathfrak{a}}_{1}(N)$ are finitely generated $\widehat{R}^{\mathfrak{a}}$-modules, the exact sequence \eqref{eq:2.1.3} shows that $H^{\mathfrak{a}}_{0}\left(\coker f \right)$ and $H^{\mathfrak{a}}_{1}\left(\coker f \right)$ are finitely generated  $\widehat{R}^{\mathfrak{a}}$-modules, and thus $\coker f$ is $\mathfrak{a}$-cofinite by (i). From the short exact sequence \eqref{eq:2.1.2}, we conclude that $\im f$ is $\mathfrak{a}$-cofinite, and from the short exact sequence \eqref{eq:2.1.1}, we infer that $\ker f$ is $\mathfrak{a}$-cofinite. It follows that $\mathcal{M}(R,\mathfrak{a})_{cof}$ is an abelian subcategory of $\mathcal{M}(R)$.
\end{prf}

It is well-known that $\cd(\mathfrak{a},R) \leq \ara(\mathfrak{a})$. On the other hand, the following example shows that an ideal $\mathfrak{a}$ of $R$ with $\cd(\mathfrak{a},R)=1$ need not have $\ara(\mathfrak{a})=1$. Hence Theorem \ref{2.1} (ii) genuinely generalizes Theorem \ref{1.4}.

\begin{example} \label{2.2}
Let $k$ be a field and $S=k[[X,Y,Z,W]]$. Consider the elements $f=XW-YZ$, $g=Y^{3}-X^{2}Z$, and $h=Z^{3}-Y^2W$ of $S$. Let $R=S/fS$, and $\mathfrak{a}=(f,g,h)S/fS$. Then $R$ is a noetherian local ring of dimension $3$, $\cd(\mathfrak{a},R)=1$, and $\ara(\mathfrak{a})\geq 2$. See \cite[Remark 2.1 (ii)]{HeSt}.
\end{example}

\section{Question 1.3}

In this section, we exploit the technique of way-out functors as the main tool to depart from modules to complexes.

\begin{definition} \label{3.1}
Let $R$ and $S$ be two rings, and $\mathcal{F}: \mathcal{D}(R) \rightarrow \mathcal{D}(S)$ a covariant functor. We say that
\begin{enumerate}
\item[(i)] $\mathcal{F}$ is \textit{way-out left} if for every $n \in \mathbb{Z}$, there is an $m \in \mathbb{Z}$, such that for any $R$-complex $X$ with $\sup X \leq m$, we have $\sup \mathcal{F}(X) \leq n$.
\item[(ii)] $\mathcal{F}$ is \textit{way-out right} if for every $n \in \mathbb{Z}$, there is an $m \in \mathbb{Z}$, such that for any $R$-complex $X$ with $\inf X \geq m$, we have $\inf \mathcal{F}(X) \geq n$.
\item[(iii)] $\mathcal{F}$ is \textit{way-out} if it is both way-out left and way-out right.
\end{enumerate}
\end{definition}

The Way-out Lemma appears in \cite[Ch. I, Proposition 7.3]{Ha2}. However, we need a refined version which is tailored to our needs. Since the proof of the original result in \cite[Ch. I, Proposition 7.3]{Ha2} is left to the reader, we deem it appropriate to include a proof of our refined version for the convenience of the reader as well as bookkeeping.

\begin{lemma} \label{3.2}
Let $R$ and $S$ be two rings, and $\mathcal{F}: \mathcal{D}(R) \rightarrow \mathcal{D}(S)$ a triangulated covariant functor. Let $\mathcal{A}$ be an additive subcategory of $\mathcal{M}(R)$, and $\mathcal{B}$ an abelian subcategory of $\mathcal{M}(S)$ which is closed under extensions. Suppose that $H_{i}\left(\mathcal{F}(M)\right) \in \mathcal{B}$ for every $M \in \mathcal{A}$ and every $i \in \mathbb{Z}$. Then the following assertions hold:
\begin{enumerate}
\item[(i)] If $X \in \mathcal{D}_{\square}(R)$ with $H_{i}(X) \in \mathcal{A}$ for every $i \in \mathbb{Z}$, then $H_{i}\left(\mathcal{F}(X)\right) \in \mathcal{B}$ for every $i \in \mathbb{Z}$.
\item[(ii)] If $\mathcal{F}$ is way-out left and $X \in \mathcal{D}_{\sqsubset}(R)$ with $H_{i}(X) \in \mathcal{A}$ for every $i \in \mathbb{Z}$, then $H_{i}\left(\mathcal{F}(X)\right) \in \mathcal{B}$ for every $i \in \mathbb{Z}$.
\item[(iii)] If $\mathcal{F}$ is way-out right and $X \in \mathcal{D}_{\sqsupset}(R)$ with $H_{i}(X) \in \mathcal{A}$ for every $i \in \mathbb{Z}$, then $H_{i}\left(\mathcal{F}(X)\right) \in \mathcal{B}$ for every $i \in \mathbb{Z}$.
\item[(iv)] If $\mathcal{F}$ is way-out and $X \in \mathcal{D}(R)$ with $H_{i}(X) \in \mathcal{A}$ for every $i \in \mathbb{Z}$, then $H_{i}\left(\mathcal{F}(X)\right) \in \mathcal{B}$ for every $i \in \mathbb{Z}$.
\end{enumerate}
\end{lemma}

\begin{prf}
(i): Let $s= \sup(X)$. Since $\amp(X) < \infty$, we argue by induction on $n= \amp(X)$. If $n=0$, then $X \simeq \Sigma^{s} H_{s}(X)$. Therefore,
$$H_{i}\left(\mathcal{F}(X)\right) \cong H_{i}\left(\mathcal{F}\left(\Sigma^{s} H_{s}(X)\right)\right) \cong H_{i-s}\left(\mathcal{F}\left(H_{s}(X)\right)\right) \in \mathcal{B},$$
as $H_{s}(X) \in \mathcal{A}$. Now, let $n \geq 1$ and assume that the result holds for amplitude less than $n$. Since $X \simeq X_{s \subset}$, there is a distinguished triangle
\begin{equation} \label{eq:3.2.1}
\Sigma^{s}H_{s}(X) \rightarrow X \rightarrow X_{s-1\subset} \rightarrow.
\end{equation}
It is clear that the two $R$-complexes $\Sigma^{s}H_{s}(X)$ and $X_{s-1\subset}$ have all their homology modules in $\mathcal{A}$
and their amplitudes are less than $n$. Therefore, the induction hypothesis implies that $H_{i}\left(\mathcal{F}\left(\Sigma^{s}H_{s}(X)\right)\right) \in \mathcal{B}$ and $H_{i}\left(\mathcal{F}\left(X_{s-1\subset}\right)\right) \in \mathcal{B}$ for every $i \in \mathbb{Z}$. Applying the functor $\mathcal{F}$ to the distinguished triangle \eqref{eq:3.2.1}, we get the distinguished triangle
$$\mathcal{F}\left(\Sigma^{s}H_{s}(X)\right) \rightarrow \mathcal{F}(X) \rightarrow \mathcal{F}(X_{s-1\subset}) \rightarrow,$$
which in turn yields the long exact homology sequence
$$\cdots \rightarrow H_{i+1}\left(\mathcal{F}(X_{s-1\subset})\right) \rightarrow H_{i}\left(\mathcal{F}\left(\Sigma^{s}H_{s}(X)\right)\right) \rightarrow H_{i}\left(\mathcal{F}(X)\right) \rightarrow $$$$ H_{i}\left(\mathcal{F}(X_{s-1\subset})\right) \rightarrow H_{i-1}\left(\mathcal{F}\left(\Sigma^{s}H_{s}(X)\right)\right) \rightarrow \cdots.$$
We break the displayed part of the above exact sequence into the following exact sequences
$$H_{i+1}\left(\mathcal{F}(X_{s-1\subset})\right) \rightarrow H_{i}\left(\mathcal{F}\left(\Sigma^{s}H_{s}(X)\right)\right)\rightarrow K\rightarrow 0,$$
$$0\rightarrow K\rightarrow H_{i}\left(\mathcal{F}(X)\right)\rightarrow L\rightarrow 0,$$
$$0\rightarrow L \rightarrow H_{i}\left(\mathcal{F}(X_{s-1\subset})\right)\rightarrow H_{i-1}\left(\mathcal{F}\left(\Sigma^{s}H_{s}(X)\right)\right).$$
Since the subcategory $\mathcal{B}$ is abelian, we conclude from the first and the third exact sequences above that $K,L \in \mathcal{B}$. Since $\mathcal{B}$ is closed under extensions, the second exact sequence above implies that $H_{i}\left(\mathcal{F}(X)\right) \in \mathcal{B}$ for every $i \in \mathbb{Z}$.

(ii): Let $i \in \mathbb{Z}$. Since $\mathcal{F}$ is way-out left, we can choose an integer $j \in \mathbb{Z}$ corresponding to $i-1$. Apply the functor $\mathcal{F}$ to the distinguished triangle
$$X_{\supset j+1} \rightarrow X \rightarrow X_{j \subset} \rightarrow,$$
to get the distinguished triangle
$$\mathcal{F}(X_{\supset j+1}) \rightarrow \mathcal{F}(X) \rightarrow \mathcal{F}(X_{j \subset}) \rightarrow.$$
From the associated long exact homology sequence, we get
$$0= H_{i+1}\left(\mathcal{F}(X_{j \subset})\right) \rightarrow H_{i}\left(\mathcal{F}(X_{\supset j+1})\right) \rightarrow H_{i}\left(\mathcal{F}(X)\right) \rightarrow H_{i}\left(\mathcal{F}(X_{j \subset})\right) =0,$$
where the vanishing is due to the choice of $j$. Since $X_{\supset j+1} \in \mathcal{D}_{\square}(R)$ with $H_{i}(X_{\supset j+1}) \in \mathcal{A}$ for every $i \in \mathbb{Z}$, it follows from (i) that $H_{i}\left(\mathcal{F}(X_{\supset j+1})\right) \in \mathcal{B}$ for every $i \in \mathbb{Z}$, and as a consequence, $H_{i}\left(\mathcal{F}(X)\right) \in \mathcal{B}$ for every $i \in \mathbb{Z}$.

(iii): Given $i \in \mathbb{Z}$, choose the integer $j$ corresponding to $i+1$. The rest of the proof is similar to (ii) using the  distinguished triangle $$X_{\supset j} \rightarrow X \rightarrow X_{j-1\subset} \rightarrow.$$

(iv): Apply the functor $\mathcal{F}$ to the distinguished triangle
$$X_{\supset 1} \rightarrow X \rightarrow X_{0 \subset} \rightarrow,$$
to get the distinguished triangle
$$\mathcal{F}(X_{\supset 1}) \rightarrow \mathcal{F}(X) \rightarrow \mathcal{F}(X_{0 \subset}) \rightarrow.$$
Since $X_{0 \subset} \in \mathcal{D}_{\sqsubset}(R)$ and $X_{\supset 1} \in \mathcal{D}_{\sqsupset}(R)$ with $H_{i}(X_{0 \subset}), H_{i}(X_{\supset 1}) \in \mathcal{A}$ for every $i \in \mathbb{Z}$, we deduce from (ii) and (iii) that $H_{i}\left(\mathcal{F}(X_{0 \subset})\right), H_{i}\left(\mathcal{F}(X_{\supset 1})\right) \in \mathcal{B}$ for every $i \in \mathbb{Z}$. Using the associated long exact homology sequence, an argument similar to (i) yields that $H_{i}\left(\mathcal{F}(X)\right) \in \mathcal{B}$ for every $i \in \mathbb{Z}$.
\end{prf}

The next result provides us with a suitable transition device from modules to complexes when dealing with cofiniteness.

\begin{theorem} \label{3.3}
If $\mathfrak{a}$ is an ideal of $R$, then the functor ${\bf R}\Gamma_{\mathfrak{a}}(-): \mathcal{D}(R) \rightarrow \mathcal{D}(R)$ is triangulated and way-out.
As a consequence, if $H^{i}_{\mathfrak{a}}(M)$ is $\mathfrak{a}$-cofinite for every finitely generated $R$-module $M$ and every $i \geq 0$, and $\mathcal{M}(R,\mathfrak{a})_{cof}$ is an abelian category, then $H^{i}_{\mathfrak{a}}(X)$ is $\mathfrak{a}$-cofinite for every $X \in \mathcal{D}^{f}(R)$ and every $i \in \mathbb{Z}$.
\end{theorem}

\begin{prf} By \cite[Corollary 3.1.4]{Li}, the functor ${\bf R}\Gamma_{\mathfrak{a}}(-): \mathcal{D}(R) \rightarrow \mathcal{D}(R)$ is triangulated and way-out. Now, let $\mathcal{A}$ be the subcategory of finitely generated $R$-modules, and let $\mathcal{B}:= \mathcal{M}(R,\mathfrak{a})_{cof}$. It can be easily seen that $\mathcal{B}$ is closed under extensions. It now follows from Lemma \ref{3.2} that $H^{i}_{\mathfrak{a}}(X) = H_{-i}\left({\bf R}\Gamma_{\mathfrak{a}}(X)\right) \in \mathcal{B}$ for every $X \in \mathcal{D}^{f}(R)$ and every $i \in \mathbb{Z}$.
\end{prf}

\begin{lemma} \label{3.4}
Suppose that $R$ admits a dualizing complex $D$, and $\mathfrak{a}$ is an ideal of $R$. Further, suppose that $H^{i}_{\mathfrak{a}}(Z)$ is $\mathfrak{a}$-cofinite for
every $Z\in \mathcal{D}_{\sqsubset}^{f}(R)$ and every $i\in \mathbb{Z}$. Let $Y\in \mathcal{D}_{\sqsupset}^{f}(R)$, and $X := {\bf R}\Hom_{R}\left(Y,{\bf R}
\Gamma_{\mathfrak{a}}(D)\right)$. Then $H_{i}(X)$ is $\mathfrak{a}$-cofinite for every $i \in \mathbb{Z}$.
\end{lemma}

\begin{prf}
Set $Z:= {\bf R}\Hom_{R}(Y,D)$. Then clearly, $Z \in \mathcal{D}_{\sqsubset}^{f}(R)$. Let $\check{C}(\underline{a})$ denote the \v{C}ech complex on a sequence of elements $\underline{a}=a_{1},...,a_{n}\in R$ that generates $\mathfrak{a}$. For any $R$-complex $W$, \cite[Proposition 3.1.2]{Li} yields that ${\bf R}\Gamma_{\mathfrak{a}}(W)\simeq \check{C}(\underline{a})\otimes_{R}^{\bf L} W$. Now, by applying the Tensor Evaluation Isomorphism, we get the following display:
\begin{equation*}
\begin{split}
X & = {\bf R}\Hom_{R}\left(Y,{\bf R}\Gamma_{\mathfrak{a}}(D)\right)\\
& \simeq \check{C}(\underline{a})\otimes_{R}^{\bf L} {\bf R}\Hom_R(Y,D)\\
& \simeq \check{C}(\underline{a})\otimes_{R}^{\bf L} Z\\
& \simeq {\bf R}\Gamma_{\mathfrak{a}}(Z).\\
\end{split}
\end{equation*}
Hence $H_{i}(X) \cong H^{-i}_{\mathfrak{a}}(Z)$ for every $i \in \mathbb{Z}$, and so and the conclusion follows.
\end{prf}

The next result answers Hartshorne's third question.

\begin{theorem} \label{3.5}
Let $\mathfrak{a}$ be an ideal of $R$ and $X\in \mathcal{D}_{\sqsubset}(R)$. Then the following assertions hold:
\begin{enumerate}
\item[(i)] If $H_i(X)$ is $\frak{a}$-cofinite for every $i\in \mathbb{Z}$, then $X$ is $\frak{a}$-cofinite.
\item[(ii)] Assume that $R$ admits a dualizing complex $D$, $\frak{a}$ is contained in the Jacobson radical of $R$, and $H^{i}_{\mathfrak{a}}(Z)$
is $\mathfrak{a}$-cofinite for every $Z\in \mathcal{D}_{\sqsubset}^{f}(R)$ and every $i\in \mathbb{Z}$. If $X$ is $\frak{a}$-cofinite in the sense
of Hartshorne, then $H_i(X)$ is $\frak{a}$-cofinite for all $i\in \mathbb{Z}$.
\end{enumerate}
\end{theorem}

\begin{prf}
(i) Suppose that $H_i(X)$ is $\frak{a}$-cofinite for all $i\in \mathbb{Z}$. The spectral sequence $$E_{p,q}^{2}=\Ext_{R}^{p}\left(R/\mathfrak{a},
H_{-q}(X)\right)\underset{p}\Rightarrow \Ext_{R}^{p+q}(R/\mathfrak{a},X)$$ from the proof of \cite[Proposition 6.2]{Ha1}, together with the assumption
that $E_{p,q}^{2}$ is finitely generated for every $p,q\in \mathbb{Z}$, conspire to imply that $\Ext_{R}^{p+q}(R/\mathfrak{a},X)$ is finitely generated.
On the other hand, one has $$\Supp_R(X)\subseteq \underset{i\in \mathbb{Z}}\bigcup \Supp_R(H_i(X))\subseteq \V(\frak a).$$ Thus $X$ is $\mathfrak{a}$-cofinite.

(ii) Suppose that $X$ is $\frak{a}$-cofinite in the sense of Hartshorne. Then by definition, there is $Y\in  \mathcal{D}^f(R)$ such that $X\simeq {\bf R}
\Hom_{R}\left(Y,{\bf R} \Gamma_{\mathfrak{a}}(D)\right)$. Now, the Affine Duality Theorem \cite[Theorem 4.3.1]{Li} implies that $$Y\otimes_{R}^{\bf L}
\widehat{R}^{\mathfrak{a}}\simeq {\bf R}\Hom_{R}\left(X,{\bf R} \Gamma_{\mathfrak{a}}(D)\right).$$ Since  $\id_R({\bf R} \Gamma_{\mathfrak{a}}(D))<\infty$
and $X\in \mathcal{D}_{\sqsubset}(R)$, we conclude that ${\bf R}\Hom_{R}\left(X,{\bf R} \Gamma_{\mathfrak{a}}(D)\right)\in \mathcal{D}_{\sqsupset}(R)$.
As the functor $-\otimes_{R}\widehat{R}^{\mathfrak{a}}:\mathcal{M}(R)\lo  \mathcal{M}(R)$ is faithfully flat, it turns out that $Y\in  \mathcal{D}_{\sqsupset}^f(R)$.
Now, the claim follows by Lemma \ref{3.4}.
\end{prf}

\begin{corollary} \label{3.6}
Let $\mathfrak{a}$ be an ideal of $R$ such that either $\cd(\fa,R)\leq 1$, or $\dim R/\fa\leq 1$, or $\dim(R)\leq 2$. Then the following assertions hold:
\begin{enumerate}
\item[(i)] $H^{i}_{\mathcal{\fa}}(X)$ is $\mathfrak{a}$-cofinite for every $X \in \mathcal{D}^{f}(R)$ and every $i \in \mathbb{Z}$.
\item[(ii)] Assume that $R$ admits a dualizing complex $D$ and $\frak{a}$ is contained in the Jacobson radical of $R$. If $X\in \mathcal{D}_{\sqsubset}(R)$
is $\frak{a}$-cofinite in the sense of Hartshorne, then $H_i(X)$ is $\frak{a}$-cofinite for all $i\in \mathbb{Z}$.
\end{enumerate}
\end{corollary}

\begin{prf}
(i) Follows from Theorem \ref{1.4}, \cite[Corollary 3.14]{Me2}, Theorem \ref{2.1} (ii) and Theorem \ref{3.3}.

(ii) Follows by (i) and Theorem \ref{3.5} (ii).
\end{prf}

\section{Correlation between Questions 1.1, 1.2 and 1.3}

In this section, we probe the connection between Hartshorne's questions as highlighted in the Introduction.

Some special cases of the following result is more or less proved in \cite[Theorem 3.10 and Proposition 3.13]{PSY}. However, we include it here with a different and
shorter proof due to its pivotal role in the theory of cofiniteness.

\begin{lemma} \label{4.1}
Let $\mathfrak{a}$ be an ideal of $R$ and $X\in \mathcal{D}_{\square}(R)$. Then the following assertions are equivalent:
\begin{enumerate}
\item[(i)] ${\bf R}\Hom_{R}(R/\mathfrak{a},X) \in \mathcal{D}^{f}(R)$.
\item[(ii)] ${\bf R}\Gamma_{\mathfrak{a}}(X)\simeq {\bf R}\Gamma_{\mathfrak{a}}(Z)$ for some $Z \in \mathcal{D}^{f}_{\square}\left(\widehat{R}^{\mathfrak{a}}\right)$.
\item[(iii)] ${\bf R}\Gamma_{\mathfrak{a}}(X)\simeq {\bf R}\Hom_{\widehat{R}^{\mathfrak{a}}}\left(Y,{\bf R}\Gamma_{\mathfrak{a}}(D)\right)$ for some
$Y \in \mathcal{D}^{f}_{\square}\left(\widehat{R}^{\mathfrak{a}}\right)$, provided that $\widehat{R}^{\mathfrak{a}}$ enjoys a dualizing complex $D$.
\end{enumerate}
\end{lemma}

\begin{prf}
(i) $\Rightarrow$ (ii): By Lemma \ref{2.0}, $Z:={\bf L}\Lambda^{\mathfrak{a}}(X) \in \mathcal{D}^{f}_{\square}\left(\widehat{R}^{\mathfrak{a}}\right)$.
Then by \cite[Corollary after (0.3)$^{\ast}$]{AJL}, we have
$${\bf R}\Gamma_{\mathfrak{a}}(Z)\simeq {\bf R}\Gamma_{\mathfrak{a}}\left({\bf L}\Lambda^{\mathfrak{a}}(X)\right)\simeq {\bf R}\Gamma_{\mathfrak{a}}(X).$$

(ii) $\Rightarrow$ (iii): Set $Y:={\bf R}\Hom_{\widehat{R}^{\mathfrak{a}}}(Z,D)$. If $\id_{\widehat{R}^{\mathfrak{a}}}(D) = n$, then there is a semi-injective resolution $D \xrightarrow {\simeq} I$ of $D$ such that $I_{i}=0$ for every $i > \sup D$ or $i < -n$. In particular, $I$ is bounded. On the other hand, $Z \in \mathcal{D}^{f}_{\square}\left(\widehat{R}^{\mathfrak{a}}\right)$, so there is a bounded $\widehat{R}^{\mathfrak{a}}$-complex $Z^{\prime}$ such that $Z\simeq Z^{\prime}$. Therefore,
$$Y= {\bf R}\Hom_{\widehat{R}^{\mathfrak{a}}}(Z,D) \simeq {\bf R}\Hom_{\widehat{R}^{\mathfrak{a}}}(Z^{\prime},D) \simeq \Hom_{\widehat{R}^{\mathfrak{a}}}(Z^{\prime},I).$$
But it is obvious that $\Hom_{\widehat{R}^{\mathfrak{a}}}(Z^{\prime},I)$ is bounded, so $Y\in \mathcal{D}^{f}_{\square}\left(\widehat{R}^{\mathfrak{a}}\right)$.

Now, let $\check{C}(\underline{a})$ denote the \v{C}ech complex on a sequence of elements $\underline{a}=a_{1},...,a_{n}\in R$ that generates $\mathfrak{a}$.
We have
\begin{equation*}
\begin{split}
{\bf R}\Gamma_{\mathfrak{a}}(X) & \simeq {\bf R}\Gamma_{\mathfrak{a}}(Z)\\
& \simeq {\bf R}\Gamma_{\mathfrak{a}}\left({\bf R}\Hom_{\widehat{R}^{\mathfrak{a}}}\left({\bf R}\Hom_{\widehat{R}^{\mathfrak{a}}}(Z,D),D\right)\right)\\
& \simeq {\bf R}\Gamma_{\mathfrak{a}}\left({\bf R}\Hom_{\widehat{R}^{\mathfrak{a}}}\left(Y,D\right)\right)\\
& \simeq \check{C}(\underline{a})\otimes_{R}^{\bf L} {\bf R}\Hom_{\widehat{R}^{\mathfrak{a}}}\left(Y,D\right)\\
& \simeq {\bf R}\Hom_{\widehat{R}^{\mathfrak{a}}}\left(Y, \check{C}(\underline{a})\otimes_{R}^{\bf L} D\right)\\
& \simeq {\bf R}\Hom_{\widehat{R}^{\mathfrak{a}}}\left(Y,{\bf R}\Gamma_{\mathfrak{a}}(D)\right).\\
\end{split}
\end{equation*}
The second isomorphism is due to the fact that $D$ is a dualizing $\widehat{R}^{\mathfrak{a}}$-module, and the fifth isomorphism follows from the application of the Tensor Evaluation Isomorphism. The other isomorphisms are straightforward.

(iii) $\Rightarrow$ (i): Similar to the argument of the implication (ii) $\Rightarrow$ (iii), we conclude that ${\bf R}\Hom_{\widehat{R}^{\mathfrak{a}}}(Y,D)\in \mathcal{D}^{f}_{\square}\left(\widehat{R}^{\mathfrak{a}}\right)$. We further have
\begin{equation*}
\begin{split}
{\bf L}\Lambda^{\mathfrak{a}}(X) & \simeq {\bf L}\Lambda^{\mathfrak{a}}\left({\bf R}\Gamma_{\mathfrak{a}}(X)\right)\\
& \simeq {\bf L}\Lambda^{\mathfrak{a}}\left({\bf R}\Hom_{\widehat{R}^{\mathfrak{a}}}\left(Y,{\bf R}\Gamma_{\mathfrak{a}}(D)\right)\right)\\
& \simeq {\bf L}\Lambda^{\mathfrak{a}}\left({\bf R}\Gamma_{\mathfrak{a}}\left({\bf R}\Hom_{\widehat{R}^{\mathfrak{a}}}(Y,D)\right)\right)\\
& \simeq {\bf L}\Lambda^{\mathfrak{a}}\left({\bf R}\Hom_{\widehat{R}^{\mathfrak{a}}}(Y,D)\right)\\
& \simeq {\bf R}\Hom_{\widehat{R}^{\mathfrak{a}}}(Y,D)\in \mathcal{D}^{f}_{\square}\left(\widehat{R}^{\mathfrak{a}}\right).\\
\end{split}
\end{equation*}
The first and the fourth isomorphisms use \cite[Corollary after (0.3)$^{\ast}$]{AJL}, the third isomorphism follows from the application of the Tensor Evaluation Isomorphism
just as in the previous paragraph, and the fifth isomorphism follows from \cite[Theorem 1.21]{PSY}, noting that as ${\bf R}\Hom_{\widehat{R}^{\mathfrak{a}}}(Y,D)\in \mathcal{D}^{f}_{\square}\left(\widehat{R}^{\mathfrak{a}}\right)$, its homology modules are $\mathfrak{a}$-adically complete $\widehat{R}^{\mathfrak{a}}$-modules. Now, the results follows from Lemma \ref{2.0}.
\end{prf}

\begin{corollary} \label{4.2}
Let $\mathfrak{a}$ be an ideal of $R$ for which $R$ is $\mathfrak{a}$-adically complete and $X\in \mathcal{D}_{\square}(R)$. Then the following assertions are equivalent:
\begin{enumerate}
\item[(i)] $X$ is $\mathfrak{a}$-cofinite.
\item[(ii)] $X\simeq {\bf R}\Gamma_{\mathfrak{a}}(Z)$ for some $Z\in \mathcal{D}^{f}_{\square}(R)$.
\item[(iii)] $X\simeq {\bf R}\Hom_{R}\left(Y,{\bf R}\Gamma_{\mathfrak{a}}(D)\right)$ for some $Y\in \mathcal{D}^{f}_{\square}(R)$, provided that $R$ enjoys a dualizing complex $D$.
\end{enumerate}
\end{corollary}

\begin{prf}
For any two $R$-complexes $V\in \mathcal{D}^{f}_{\square}(R)$ and $W\in \mathcal{D}_{\square}(R)$, one may easily see that $$\Supp_{R}
\left({\bf R}\Hom_{R}\left(V,{\bf R}\Gamma_{\mathfrak{a}}(W)\right)\right) \subseteq V(\mathfrak{a}).$$ Also, for any $U\in \mathcal{D}(R)$,
\cite[Corollary 3.2.1]{Li} yields that $\Supp_{R}(U)\subseteq V(\mathfrak{a})$ if and only if ${\bf R}\Gamma_{\mathfrak{a}}(U)\simeq U$.
Hence the assertions follow from Lemma \ref{4.1}.
\end{prf}

The next result reveals the correlation between Hartshorne's questions.

\begin{theorem} \label{4.3}
Let $\mathfrak{a}$ be an ideal of $R$. Consider the following assertions:
\begin{enumerate}
\item[(i)] $H^{i}_{\mathfrak{a}}(M)$ is $\mathfrak{a}$-cofinite for every finitely generated $R$-module $M$ and every $i \geq 0$, and $\mathcal{M}(R,\mathfrak{a})_{cof}$ is an abelian subcategory of $\mathcal{M}(R)$.
\item[(ii)] $H^{i}_{\mathfrak{a}}(X)$ is $\mathfrak{a}$-cofinite for every $X\in \mathcal{D}^{f}(R)$ and every $i \in \mathbb{Z}$.
\item[(iii)] An $R$-complex $X\in \mathcal{D}_{\square}(R)$ is $\mathfrak{a}$-cofinite if and only if $H_{i}(X)$ is $\mathfrak{a}$-cofinite for every $i \in \mathbb{Z}$.
\end{enumerate}
Then the implications $(i)\Rightarrow (ii)$ and $(iii)\Rightarrow (i)$ hold. Furthermore, if $R$ is $\mathfrak{a}$-adically complete, then all three assertions are equivalent.
\end{theorem}

\begin{prf}
(i) $\Rightarrow$ (ii): Follows from Theorem \ref{3.3}.

(iii) $\Rightarrow$ (i): Let $M$ be a finitely generated $R$-module. Since $H^{i}_{\mathfrak{a}}(M)=0$ for every $i<0$ or $i> \ara(\mathfrak{a})$, we have
${\bf R}\Gamma_{\mathfrak{a}}(M)\in \mathcal{D}_{\square}(R)$. However, \cite[Proposition 3.2.2]{Li} implies that
$${\bf R}\Hom_{R}\left(R/\mathfrak{a},{\bf R}\Gamma_{\mathfrak{a}}(M)\right) \simeq {\bf R}\Hom_{R}(R/\mathfrak{a},M),$$
showing that ${\bf R}\Gamma_{\mathfrak{a}}(M)$ is $\mathfrak{a}$-cofinite. The hypothesis now implies that $H^{i}_{\mathfrak{a}}(M)=H_{-i}({\bf R}\Gamma_{\mathfrak{a}}\left(M)\right)$ is $\mathfrak{a}$-cofinite for every $i\geq 0$.

Now, let $M$ and $N$ be two $\mathfrak{a}$-cofinite $R$-modules and $f:M\rightarrow N$ an $R$-homomorphism. Let $\varphi:M \rightarrow N$ be the morphism in $\mathcal{D}(R)$ represented by the roof diagram $M \xleftarrow{1^{M}} M \xrightarrow{f} N$. From the long exact homology sequence associated to the distinguished triangle
\begin{equation} \label{eq:4.3.1}
M \xrightarrow {\varphi} N \rightarrow \Cone(f) \rightarrow,
\end{equation}
we deduce that $\Supp_{R}\left(\Cone(f)\right) \subseteq \V(\mathfrak{a})$.
In addition, applying the functor ${\bf R}\Hom_{R}(R/\mathfrak{a},-)$ to \eqref{eq:4.3.1}, gives the distinguished triangle
$${\bf R}\Hom_{R}(R/\mathfrak{a},M) \rightarrow {\bf R}\Hom_{R}(R/\mathfrak{a},N) \rightarrow {\bf R}\Hom_{R}\left(R/\mathfrak{a},\Cone(f)\right)
\rightarrow,$$
whose associated long exact homology sequence shows that
$${\bf R}\Hom_{R}\left(R/\mathfrak{a},\Cone(f)\right)\in \mathcal{D}^{f}(R).$$
Hence, the $R$-complex $\Cone(f)$ is $\mathfrak{a}$-cofinite. However, we have
$$\Cone(f)= \cdots \rightarrow 0 \rightarrow M \xrightarrow {f} N \rightarrow 0 \rightarrow \cdots,$$
so $\Cone(f)\in \mathcal{D}_{\square}(R)$. Thus the hypothesis implies that $H_{i}\left(\Cone(f)\right)$ is $\mathfrak{a}$-cofinite for every
$i\in \mathbb{Z}$. It follows that $\ker f$ and $\coker f$ are $\mathfrak{a}$-cofinite, and as a consequence $\mathcal{M}(R,\mathfrak{a})_{cof}$ is
an abelian subcategory of $\mathcal{M}(R)$.

Now, suppose that $R$ is $\mathfrak{a}$-adically complete.

(ii) $\Rightarrow$ (iii): Let $X\in \mathcal{D}_{\square}(R)$. Suppose that $H_{i}(X)$ is $\mathfrak{a}$-cofinite for every $i \in \mathbb{Z}$.
Then Theorem \ref{3.5} (i) yields that $X$ is $\mathfrak{a}$-cofinite.

Conversely, assume that $X$ is $\mathfrak{a}$-cofinite. Then by Corollary
\ref{4.2}, $X\simeq {\bf R}\Gamma_{\mathfrak{a}}(Z)$ for some $Z\in \mathcal{D}^{f}_{\square}(R)$. Thus the hypothesis implies that $$H_{i}(X)
\cong H_{i}\left({\bf R}\Gamma_{\mathfrak{a}}(Z)\right)=H^{-i}_{\mathfrak{a}}(Z)$$ is $\mathfrak{a}$-cofinite for every $i \in \mathbb{Z}$.
\end{prf}

In view of Corollary \ref{4.2}, the next result answers Hartshorne's third question for homologically bounded $R$-complexes.

\begin{corollary} \label{4.4}
Let $\mathfrak{a}$ be an ideal of $R$ for which $R$ is $\mathfrak{a}$-adically complete. Suppose that either $\cd(\mathfrak{a},R) \leq 1$, or
$\dim(R/\mathfrak{a}) \leq 1$, or $\dim(R) \leq 2$. Then an $R$-complex $X\in \mathcal{D}_{\square}(R)$ is $\mathfrak{a}$-cofinite if and only
if $H_{i}(X)$ is $\mathfrak{a}$-cofinite for every $i \in \mathbb{Z}$.
\end{corollary}

\begin{prf}
Obvious in light of Corollary \ref{3.6} and Theorem \ref{4.3}.
\end{prf}

\begin{acknowledgement}
The authors are deeply grateful to Professor Robin Hartshorne for his invaluable comments on an earlier draft of this paper.
\end{acknowledgement}

%%%%%%%%%%%%%%%%%%%%%%%%%%%%%%%%%%%%%%%%%%%%%%%%%%%%%%%%%%%%%%%%%%%%%%%%%%%%%%%%%%%%%%%%%%%%%%%%%%%%%%%%%%%%%%%%%%%%%%%%%%%%%%%%%%%%%%%%%%%%%%%%%%%%%%%%%%%%%%%%%


\begin{thebibliography}{99}

\bibitem[AJL]{AJL}{L. Alonso Tarr\'{i}o, A. Jerem\'{i}as L\'{o}pez and J. Lipman}, {\it Local homology and cohomology on
schemes}, Ann. Sci. \'{E}cole Norm. Sup., (4), \textbf{30}(1), (1997), 1-39.

\bibitem[AF]{AF}{L. Avramov and H-B. Foxby}, {\it Homological dimensions of unbounded complexes}, J. Pure Appl.
Algebra, \textbf{71}(2-3), (1991), 129-155.

\bibitem[BNS]{BNS}{K. Bahmanpour, R. Naghipour and M. Sedghi}, {\it On the category of cofinite modules which is abelian},
Proc. Amer. Math. Soc., \textbf{142}(4), (2014), 1101-1107.

\bibitem[BN]{BN}{K. Bahmanpour and R. Naghipour}, {\it Cofiniteness of local cohomology modules for ideals of small
dimension}, J. Algebra, \textbf{321}(7), (2009), 1997-2011.

\bibitem[DM]{DM}{D. Delfino and T. Marley}, {\it Cofinite modules and local cohomology}, J. Pure Appl. Algebra,
\textbf{121}(1), (1997), 45-52.

\bibitem[EK]{EK}{K. Eto and K. Kawasaki}, {\it A characterization of cofinite complexes over complete Gorenstein domains},
J. Commut. Algebra, \textbf{3}(4), (2011),  537-550.

\bibitem[Fo]{Fo}{H-B. Foxby}, {\it Hyperhomological algebra and commutative rings}, notes in preparation.

\bibitem[GM]{GM}{J. P. C. Greenlees and J. P. May}, {\it Derived functors of $I$-adic completion and local homology},
J. Algebra, \textbf{149}(2), (1992), 438-453.

\bibitem[Ha1]{Ha1}{R. Hartshorne}, {\it Affine duality and cofiniteness}, Invent. Math., \textbf{9}, (1969/1970), 145-164.

\bibitem[Ha2]{Ha2}{R. Hartshorne}, {\it Residues and duality}, Lecture notes in mathematics, \textbf{20}, (1966).

\bibitem[HeSt]{HeSt}{M. Hellus and J. Str\"{u}kard}, {\it Local cohomology and Matlis duality}, Univ. Iagel. Acta Math., \textbf{45},
(2007), 63-70.

\bibitem[HK]{HK}{G. Huneke and J. Koh}, {\it Cofiniteness and vanishing of local cohomology modules}, Math. Proc.
Camb. Phil. Soc., \textbf{110}(3), (1991), 421-429.

\bibitem[HS]{HS}{G. Huneke and R.Y.  Sharp}, {\it Bass numbers of local cohomology modules}, Trans. Amer. Math. Soc.,
\textbf{339}(2), (1993), 765-779.

\bibitem[Ka1]{Ka1}{K. Kawasaki}, {\it On the category of cofinite modules for principal ideals}, Nihonkai Math. J.,
\textbf{22}(2), (2011), 67-71.

\bibitem[Ka2]{Ka2}{K. Kawasaki}, {\it Cofiniteness of local cohomology modules for principal ideals}, Bull. London
Math. Soc., \textbf{30}(3), (1998), 241-246.

\bibitem[Li]{Li}{J. Lipman}, {\it Lectures on local cohomology and duality}, Local cohomology and its applications,
Lecture notes in pure and applied mathematics, \textbf{226}, (2012), Marcel Dekker, Inc.

\bibitem[Ly]{Ly}{G. Lyubeznik}, {\it Finiteness properties of local cohomology modules (an application of D-modules
to commutative algebra)}, Invent. Math., \textbf{113}(1), (1993), 41-55.

\bibitem[Me1]{Me1}{L. Melkersson}, {\it Cofiniteness with respect to ideals of dimension one}, J. Algebra,
\textbf{372}, (2012), 459-462.

\bibitem[Me2]{Me2}{L. Melkersson}, {\it Modules cofinite with respect to an ideal}, J. Algebra, \textbf{285}(2),
(2005), 649-668.

\bibitem[PAB]{PAB}{G. Pirmohammadi, K. Ahmadi Amoli and K. Bahmanpour}, {\it Some homological properties of ideals
with cohomological dimension one}, Colloq. Math., \textbf{149}(2), (2017), 225-238.

\bibitem[PSY]{PSY}{M. Porta, L. Shaul and A. Yekutieli}, {\it Cohomologically cofinite complexes}, Comm. Algebra,
\textbf{43}(2), (2015), 597-615.

\bibitem[Sp]{Sp}{N. Spaltenstein}, {\it Resolutions of unbounded complexes}, Compositio Math., \textbf{65}(2),
(1988), 121-154.

\bibitem[WW]{WW}{S. Sather-Wagstaff and R. Wicklein}, {\it Support and adic finiteness for complexes}, Comm. Algebra,
\textbf{45}(6), (2017), 2569-2592.

\bibitem[Y]{Y}{K. Yoshida}, {\it Cofiniteness of local cohomology modules for ideals of dimension one}, Nagoya Math. J.,
\textbf{\bf 147}, (1997), 179-191.

\end{thebibliography}
\end{document}